\documentclass[two column, 10pt]{article}
\usepackage{authblk}
\title{Optimization for Secure and Humane Border Operations}
\author[1]{Fatemeh Farajzadeh}
\author[1,2]{Andrew C. Trapp}
\affil[1]{Data Science Program, Worcester Polytechnic Institute, Worcester, MA 01609, USA}
\affil[2]{WPI Business School, Worcester Polytechnic Institute, Worcester, MA 01609, USA}
\date{}                     
\usepackage{amsmath}
\usepackage{amsfonts}
\usepackage[colorinlistoftodos]{todonotes}
\usepackage{xcolor}
\usepackage[top=1in,bottom=1in,left=1in,right=1in]{geometry}
\usepackage{multicol}

\usepackage[numbers,sort&compress]{natbib}

\usepackage{setspace}
\usepackage[none]{hyphenat}
\usepackage[hyphens]{url}  
\usepackage[pdftex,bookmarksnumbered,hidelinks,breaklinks]{hyperref}

\begin{document}
\setcounter{page}{1}
\maketitle

In 2020 tens of millions of migrants worldwide sought safety due to various push-pull factors such as economic, geographic, demographic, and societal challenges~\cite{unhcrglobaltrends2020,IOM}.
Large movements of human flows toward national borders constitute a complex humanitarian crisis that requires appropriate preparations of border resources. 
Motivated by prevailing national agenda issues for improved border support, we identify common operational challenges along the border and discuss how they can be addressed using advanced analytics.
Multiple stakeholders, corresponding goals and complex decisions call for advanced optimization techniques to holistically treat such challenges.
This chapter highlights how optimization can be employed to improve both the security and humanitarian aspects of border operations in the context of migration flows.

\section*{Introduction}
Shifting and evolving migration patterns call for border plans that feature proactive mechanisms to improve operations in more effectively using limited resources.
The European Border and Coast Guard Agency, also known as Frontex, reported the need to deploy additional resources and increase operations in the Mediterranean for three primary efforts~\cite{EU1}.
The US Customs and Border Protection (CBP) agency has reported that movements of large migrant groups have resulted in excessive use of border capacities and sacrificing the border security mission~\cite{cbp1}.
While border operations must necessarily consider security and orderliness, a large majority of countries worldwide (over 150) are signatories on the 1951 Convention Relating to the Status of Refugees~\cite{unconv51}, as well as its 1967 Protocol~\cite{UN11}.
Thus, the proper management of such operational resources involves dual aspects of security and humanity, and offers opportunities for advanced analytics to inform critical border decisions.

Researchers in the early 2000's primarily emphasized issues pertaining to homeland security operations~\cite{wright}. 
Emergency preparedness and response, border and transportation security, threat analysis, and assessment were collectively positioned as the opportunities with the greatest potential for improvement via optimization techniques. 
The use of discrete optimization to improve homeland security and disaster management has also been discussed~\cite{McLay}. 
This study classified emergency and security operations into four pre-defined components of a disaster life cycle, including mitigation, preparedness, response, and recovery.  
Moreover, basic optimization models were presented according to the disaster life cycle sequence: checkpoint screening of arrivals and network fortification, prepositioning of supplies and evacuation planning, rescue and aid distribution, and infrastructure restoration.
These efforts support at-risk migrants through securing EU borders by the deployment of high-tech surveillance systems, disrupting human smuggling and trafficking networks by effective patrolling strategies, and humanitarian relief and processing~\cite{b1,EU3}.
With the fluctuating need to accommodate waves of displaced peoples, additional facilities are needed along borders where the timely supply and delivery of essential products and services are key decisions~\cite{cbp2}.

The management of human flows at national borders constitutes a complex emergency that requires a time- and cost-efficient, fair and secure planning approach.
The planning characteristics of complex emergencies differ from natural disaster response operations that typically are sudden-onset relief efforts and involve heightened vulnerabilities that, with time, tend to subside. 
Rather, border operations in relation to migration flow management is an ongoing process, and it takes place in an ad hoc manner across multiple agencies. 
Real-time planning decisions are affected by the evolution of stochastic flows over time, including high levels of uncertainty on human needs, variability of supply cycle duration, and availability of accurate information regarding the engagement of abusive groups such as smugglers and human traffickers.  
To date, relatively few studies have leveraged analytical-based methods to improve border operations in response to shifting human flows.

The aforementioned considerations motivate new structural components and modeling choices for the design of appropriate mathematical models.
Fairness considerations such as deprivation costs must be balanced alongside traditional objectives like cost, time, and distance~\cite{HOLGUI}.
Other important modeling components include the presence of multiple decision-makers that influence any one party's control level in making unilateral decisions, as well as a high level of uncertainty and unreliability associated with planning parameters.  
Mathematical optimization appears to be the most widely used methodology in managing relief resources while addressing the aforementioned challenges under a dynamic planning environment~\cite{Gina}.
And, it can be adapted to address the problems arising from global immigration crises with the goal of properly allocating security and support resources. 
This chapter proceeds to discuss border operations in two areas, the staging of security and shelter-based resources through location analysis; and search and rescue challenges in borderlands both for saving migrants and protecting against bad actors.

\section*{Location and Allocation Challenges in Border Operations}
Decisions on locating facilities have perhaps the greatest impact on the performance of relief operations in humanitarian crises, due to the direct relationship with response times and distribution costs.
Key operational issues have been studied by researchers and linked to the broader operations research literature by referring to well-known location decision problems in collaboration with the Turkish Red Crescent~\cite{Karsu}.
A recent study thoroughly summarized various types of facility location challenges according to aspects such as the location space (continuous, network-based, discrete), the type of facilities (suppliers, distribution centers, shelters), the critical decision-makings criteria (equity, cost, reliability), and source of uncertainty (supply, demand, network)~\cite{Zehranaz}.
Although facility location problems have theoretically been addressed in isolation, in practice, location decisions and positioning challenges are coupled with allocation, inventory, distribution, and routing decisions in disaster relief operation settings~\cite{Zehranaz}.
As it relates to accommodating and processing displaced populations at international borders, scholars must address the same positioning challenges while accounting for the specific aspects of the associated intricacies and corresponding planning complexities.  

The major services of border control agencies, whenever human flows arrive, are administrative and essential aid arrangements that require three primary types of decisions.
These decisions include the location of service bases in harmony with human flow directions and volumes, capacity and inventory levels according to functions of the respective facilities, and timely service and aid distribution.
The aim of establishing borderland service bases is to focus on essential operations: safety and security, health and sanitation, temporary settlements, and processing activities for identification and registration.
The aforementioned operations will provide decision-making authorities with the basic steps of strategic and tactical actions when responding to the humanitarian crises at the borders.
 
A limited number of studies have taken a similar optimization perspective in refugee camp management and administration, seeking a cost-effective and humane planning approach.
To optimally allocate aid to refugee camps under demand and replenishment cycle uncertainty, an integer optimization model with the objective of minimizing the expected total cost of holding excess inventory, referring external refugees elsewhere, and depriving internal refugees of critical aid is proposed~\cite{az}.
Another study investigated the distribution of cash and e-vouchers to refugees as well as the routing problem of trucks providing child-friendly spaces to vulnerable refugee children. 
The authors presented a bi-objective model based on covering vehicle routing problems with integrated tours where the aim is to minimize the unsatisfied demand and traveled distance~\cite{Th1}. 
A multi-objective mixed-integer program was developed to optimally choose the locations and capacity level of short-term refugee camps in Syria to minimize maintenance and startup costs and ensure camp safety~\cite{Vatasoiu}.
A combination of machine learning and integer optimization techniques are used to assign refugees to initial locations with the aim of improving integration for both refugees and destinations~\cite{Ah}.


While border operations involve issues related to flow management, they also involve security operations to protect the border from bad actors. Not surprisingly, these two missions are intertwined.  
At the same time, in the era of emerging technologies such as surveillance towers, underground sensing technologies, drones, and other AI-based solutions, border control agencies are increasingly transforming their operations through investing in the deployment of smart tools.
From the security perspective, the use of technology enhances situational awareness, agent effectiveness, and safety by supporting border patrol agents to cover more territory for the purposes of interdiction.
From the humanitarian perspective, these technologies help agents efficiently discriminate items of interest, protect vulnerable people from criminal networks, and proactively respond to detected events.
Location optimization models can improve the process of designing an efficient surveillance system to serve both purposes. This is a relatively new area of research that has recently attracted the attention of operations researchers.  
A maximum coverage optimization model was developed to determine hubs for unmanned aerial vehicles via two subproblems, namely, the selection of hub locations and demand point allocation~\cite{Mehmet}.
Similarly, a multi-objective bi-level optimization model was developed to address the problem of sensor relocation for the purpose of maximizing the minimal expected exposure of intruders when traversing a defended border region~\cite{Aa}.
Allocation optimization problems over a rolling horizon are also of great importance in this highly variable planning setting.  
These sets of optimization problems fall under dynamic resource allocation problems that are often modeled as dynamic assignment and stochastic knapsack problems~\cite{Kleywegt}.

\section*{Search and Rescue Challenges in Border Operations}
Official reports have documented the risks of illegal border crossings in remote areas arising from a lack of processing at lawful ports of entry~\cite{cbp3,UNH}. 
Migrants on the move face substandard living conditions and many types of risks, including getting lost, dangerous terrain, natural hazards, and man-made dangers such as trafficking, fraud, and smuggling.
Border agents are responsible for patrolling the assigned borderlands searching for any potential hazards, notably, people in need of humanitarian support. 
Search and rescue practices must be agile and able to track autonomous behaviors and movements in a situation the clock is ticking.  
These practices need skilled individuals with different capabilities to operate as a multi-talented team. 
In addition, a disaster-impacted zone with a time-dependent demand rate urges rescuers to be flexible in dealing with changing circumstances. 
There exist specialized units that are trained primarily to assist injured or stranded migrants, such as Border Patrol Search, Trauma, and Rescue Unit (BORSTAR) in the United States and and Search and Rescue (SAR) in Europe~\cite{cbp4,frontex}.
Thus, the use of optimization for team planning and scheduling of various rescue missions is of particular importance. 
Although rescue operations have already risen to prominence in disaster emergency management, a significant gap exists to address such problems using optimization techniques.

Search and rescue strategies must adapt to newly revealed information of team efforts, changing circumstances, and the progress of ongoing rescue operation. 
Beyond sharing prompt information of ongoing operations updated by central planners, real-life constraints must be reflected in rescue teams tours planning, including hierarchical skill levels, collaborative rescue service, time-dependent service duration, and idling rest periods~\cite{WEX}. 
In addition to operational constraints considered in rescue operations of natural disasters, border operations management is characterized by dynamics associated with various uncertain service requests or unforeseen events that typically took place in close relation with border security operations. 
Such problem dynamics may be tackled through myopic approaches, by look-ahead procedures that rely only on available information at the time of decision-making or perhaps take into account probabilistic information concerning future events.    

Search and rescue tasks can be classified by the type of terrain, such as air-sea rescue, urban rescue, hard-to-access rescue such as mountainous terrain, and general ground rescue.
Borderland search and rescue may feature additional unique characteristics, such as incident type, level of complexity, geographical factors, and immediate survival resources, all in a traumatic situation.    
Respecting geographical factors, the decision of rescue team location in sea rescue operations differs based upon the given required operation capacity and boat workload~\cite{Razi}. 
This study developed an incident-based allocation optimization model to assign responsibility zones of the rescue team considering historical demand locations to account for demand uncertainty.
A relocation possibilities of resources, using a scenario-based multi-objective programming model to generate a balanced workload distribution in the scheme of location-allocation analysis also has been addressed~\cite{Karatas1}.
Another study used a two-stage stochastic programming framework to fairly cover the demand of affected zones for the purpose of maximizing the utility of deployed rescue teams while also considering the transfer of resources to other affected areas or discharging decisions so as minimize the completion times of search and rescue teams~\cite{ahmadi}.
The challenges of determining an optimal strategy for dispatching rescue teams to post-disaster sites was also another subject of interest.
Relatively, a two-stage stochastic programming formulated to first seek a set of tours while maximizing the total expected rescued people and subsequently, second-stage variables update the allocation decisions over time by the accessibility of newly released resources and revealed planning information~\cite{CHEN}.
To enhance the performance of rescue operations, multimodal interaction is presented as an approach to incorporate a mixed-use. 
A goal programming approach is proposed to dynamically optimize the location and allocation of rescuers, boats, and helicopters considering the seasonality of demand~\cite{KARATAS2}. 
The author attempted to achieve a balanced workload for all boats and helicopters such that the sum of shortage and excess operation hours are minimized via deviations from the critical incidents response time and budget limitation. 
An integrated platform of ground and aerial vehicles is also introduced to support search and rescue activities~\cite{ca}. 
\newline
The widespread implementation of cutting-edge technologies for surveillance and monitoring makes it possible to quickly assess key planning parameters, follow changes in these conditions as they evolve, and identify at-risk persons that may need rescuing.
A variety of search techniques exist to effectively cover the target area. 
These algorithms include Bayesian approaches for tracking probable target locations, Markov decision processes for unobserved underlying states, and greedy heuristics where minimizing the solving time is the optimization criterion~\cite{waharte}.    
A study on creating a \textit{good} drone presented various successes and practical factors adapted to more effectively detect migrants, which can help researchers familiarize themselves with the operational challenges of drones~\cite{b2}.
A study explored simulation-optimization techniques to address the potential use of drones in searching for and locating victims of search and rescue providers in a mountain environment~\cite{KARACA}.
\newline
We believe this chapter outlines key borderland operations challenges in location, allocation, and search and rescue. Guided by dual humanitarian and security purposes, in conjunction with increasing technological capabilities, we sincerely hope operations research and analytics researchers consider the excellent opportunity to further explore this complex planning environment and underlying 
optimization challenges.


\section*{Acknowledgements}
We are grateful for the assistance of many in making this chapter successful. In particular, the authors are grateful for the support of the National Science Foundation (Operations Engineering) grant CMMI-1825348.


\bibliographystyle{unsrtnat} 
\bibliography{lit}

\end{document}